\documentclass[a4paper, 11pt]{article}

\usepackage{amsmath}
\usepackage{amssymb}
\usepackage{graphics}
\usepackage[dvips]{graphicx}
\usepackage{epsfig}
\usepackage[english]{babel}
\usepackage{lineno}
\usepackage{color}
\usepackage{caption}
\usepackage{subcaption}
\usepackage{epstopdf}
\graphicspath{ {./figures/} }

\newcommand{\be}{\begin{equation}}
\newcommand{\ee}{\end{equation}}
\newcommand{\bea}{\begin{eqnarray}}
\newcommand{\eea}{\end{eqnarray}}

\newcommand{\nod}{\noindent}

\newcommand{\ba}{\begin{array}}
\newcommand{\ea}{\end{array}}
\newcommand{\bc}{\begin{center}}
\newcommand{\ec}{\end{center}}

\DeclareMathOperator{\Tr}{Tr}
\DeclareMathOperator{\sign}{sign}

\newtheorem{defi}{Definition}

\begin{document}

\title{Dynamic control of modern, network-based epidemic models}

\author{Fanni S\'{e}lley$^{1}$, \'{A}d\'{a}m Besenyei$^{1}$, Istvan Z. Kiss $^{2}$ \& P\'{e}ter L. Simon $^{1,\ast}$}

\maketitle

\begin{center}
$^1$ Institute of Mathematics, E\"otv\"os Lor\'{a}nd University Budapest, and \\ Numerical Analysis and Large Networks Research Group, Hungarian Academy of Sciences, Hungary\\
$^2$ School of Mathematical and Physical Sciences, Department of Mathematics, University of Sussex, Falmer, Brighton BN1 9QH, UK
\end{center}

\vspace{1cm}

\begin{abstract}
In this paper we make the first steps to bridge the gap between classic control theory and modern,
network-based epidemic models. In particular, we apply nonlinear model predictive control (NMPC)
to a pairwise ODE model which we use to model a susceptible-infectious-susceptible ($SIS$)
epidemic on non-trivial contact structures. While classic control of epidemics concentrates on
aspects such as vaccination, quarantine and fast diagnosis, our novel setup allows us to deliver
control by altering the contact network within the population. Moreover, the ideal outcome of
control is to eradicate the disease while keeping the network well connected. The paper gives a
thorough and detailed numerical investigation of the impact and interaction of system and control
parameters on the controllability of the system. The analysis reveals, that for certain set
parameters it is possible to identify critical control bounds above which the system is
controllable. We foresee, that our approach can be extended to even more realistic or
simulation-based models with the aim to apply these to real-world situations.
\end{abstract}

\nod {\bf Keywords:} SIS epidemic; pairwise model; adaptive network; nonlinear model predictive control

\vspace{1cm}

\nod {\bf AMS classification:} 34H20, 05C82, 37N25, 92D30

\vspace{1cm}
\begin{flushleft}
$\ast$ corresponding author\\
email: simonp@cs.elte.hu\\
\end{flushleft}

\newpage

\section{Introduction}

\subsection{Background}
Being able to control a process or a systems can prove to be highly beneficial as it allows the
user to tune it or operate it in a planned or ideal regime \cite{RFC_HZ1995, EDS_1998}. Hence,
\textit{control theory} is a subject area on its own at the interface of subjects ranging from
engineering and mathematics to biology \cite{RFC_HZ1995,Grune,EDS_1998}. Mathematical models of
disease transmission, be it of simple compartmental type \cite{RMA_RMM_1991} or more modern
network-type models \cite{Danon,KeelingEames}, have been and are being developed with the ultimate
aim of making predictions about our capability to \textit{control} outbreaks. An epidemiological
model, describing disease transmission within a population, that is correctly developed and
parametrised, offers important insight  into understanding which control mechanisms and under what
circumstances can lead to a reduction in the prevalence of infection or its complete eradication.
For many models, this problem is well understood, especially in terms of vaccination
\cite{RMA_RMM_1991}, quarantine and contact tracing \cite{Kiss_CT}. However, in all this cases
control does not form an integral part of the disease dynamics and often only comes in as the
proportion of the population that needs to be vaccinated in order to develop herd immunity so that
infection can be stopped or as the critical contact tracing rate, which for sexually transmitted
infections (STIs) can differentiate between the disease-free and epidemic state.

Control, in the general sense, is dynamic in nature, where via an external input or perturbation to the system, the users are able
to tune it towards a desired outcome. This process, in many cases, is dynamic where the challenge is to determine the optimal external input across time
in order to reach a target or to minimise a cost function. In terms of epidemics, such questions have been investigated in order to determine for example the optimal time dependent vaccination
in a susceptible-infected/infectious-recovered (SIR) model under minimising a cost function that measures the cumulative amount of infected and vaccinated people \cite{MR_WKH_74}.
More recently, but still in the context of classic compartmental models, Hansen \& Day \cite{EH_TD_11} have considered optimal control in the presence of limited resources.

It is now evident that modern epidemiological models are amenable to account for and incorporate
network structure which aims to mimic to some degree a more realistic contact pattern amongst
members of a population. Pairwise models proved to be quite successful in this modelling endeavour
as they provide a relatively simple representation of epidemics unfolding on a network as opposed
to the homogenous random mixing assumption of the classic compartmental models.
 In this paper we wish to bridge the gap between modern disease transmission models \cite{HouseUnifying,KeelingEames} and control of epidemics, where the focus is on controlling the
network and not so much disease parameters, such as recovery time or the widely used pre-emptive
or reactive vaccination. This opportunity to broaden the control's target arises naturally since
the network of contact is explicitly modeled, and thus controllable. For example, in
\cite{Barabasi}, Barab\'{a}si et al. studied the controllability of complex directed networks. For a
deterministic, but not a stochastic epidemic model, they investigated how the structure of the
network influences its controllability. Their aim was to identify special vertices in the network,
the so-called driver nodes, such that the system can completely be controlled through these nodes.
By controllability the authors referred to structural controllability, which means that the system
can be controlled for almost all control values. This is a generic property of the network which
can be rephrased in terms of graph theory. By these tools the authors developed a method to find
the minimal number of driver nodes in directed networks. Then this method was applied to real
networks to study how the degree distribution of the network determines the minimal number of
driver nodes.

\subsection{The problem}
The aim of the paper is to investigate the control of a susceptible-infected/infectious-susceptible (SIS) epidemic on a network, using pairwise equations, by controlling
the creation and deletion of edges of certain types. The classic pairwise model augmented with the control elements leads to the following system of equations:
\begin{subequations}
\begin{align}
\dot{[I]}&=\tau [SI] - \gamma [I],  \label{mastereq0} \\
\dot{[SI]}&=\gamma ( [II] - [SI])+ \tau ([SSI]-[ISI]-[SI])-u_1\cdot f_1([SI]), \\
\dot{[II]}&=-2\gamma [II] +2\tau ([ISI]+[SI]), \\
\dot{[SS]}&=2\gamma [SI] - 2\tau [SSI]+u_2\cdot f_2([S],[SS]), \label{mastereq1}
\end{align}
\end{subequations}
where the $[\cdot]$ brackets denote expected number of singles and pairs of different types. For
example, $[SI]$ denotes the expected value of the number of $SI$ edges, which amounts effectively
to counting on labeled networks. The evolution equations follow naturally by observing that
singles depend on pairs, and pairs depend on triples. The precise derivation of these equations is
discussed in detail in \cite{Keeling1999}. In the system above, the control parameters are
$u_1$-the rate of cutting $SI$ edges and $u_2$-the rate of creating/deleting $SS$ edges. The
function $f_1$ and $f_2$ will be specified later, but in general these will be linear or quadratic
functions describing the precise rewiring mechanisms. The parameter $\tau$ is the per contact
infection rate and $\gamma$ is the rate of recovery. The desired outcome of our control problem is
to eradicate the epidemic, while keeping the network well connected, i.e. drive the system to
$I(T)=0$, $n(T)=n_0$, for some final time $T>0$, where $N$ is the population size and
$n(t)=([SS]+2[SI]+[II])/N$ is the average connectivity in the network.

\subsection{The structure of the paper}
The paper is organised as follows. First we consider in detail the problem of \textit{constant control}, where the problem is effectively equivalent to a dynamic or adaptive network problem, where
the epidemic dynamics and the dynamics of the network impact on and influence each other. Here, we will provide a classic bifurcation type analysis and we show that there are three qualitatively
different regimes: (1) disease-free steady state is stable, (2) stable endemic state, and finally (3) stable oscillations in both epidemic dynamics and network's average connectivity. This is followed
by the \textit{dynamic control} case, where we use the Nonlinear Model Predictive Control Method to determine if controllability is possible and how successful control depends on parameters,
such as infection rate, control bounds, the frequency of intervention and damping parameters in the control's target function. In many cases we give a substantial treatment and identify controllable
and uncontrollable situations. Finally, we discuss links to classic control and outlook towards the problem of controlling individual-based network models.

\section{Constant control}

In this section we make an attempt to control the epidemic by finding suitable values for $u_1$ and $u_2$ which stay constant until the end of the control period. We consider positive values for
these parameters, so the control removes $SI$ edges while creating new $SS$ edges. The control should delete no more edges than the existing $SI$ edges, so we take $f_1([SI])=[SI]$, and the
control should make no more $SS$ connections than the total number of unconnected $S-S$ pairs, so we take $f_2([S],[SS])=[S]([S]-1)-[SS]$.
By substituting $[S]=N-[I]$ system \eqref{mastereq0}-\eqref{mastereq1} takes the following form:
\begin{subequations}
\begin{align}
\dot{[I]}&=\tau [SI] - \gamma [I], \\
\dot{[SI]}&=\gamma ( [II] - [SI])+ \tau ([SSI]-[ISI]-[SI])-u_1[SI], \\
\dot{[II]}&=-2\gamma [II] +2\tau ([ISI]+[SI]), \\
\dot{[SS]}&=2\gamma [SI] - 2\tau [SSI]+u_2((N-[I])(N-[I]-1)-[SS]).
\end{align}
\end{subequations}
Now instead of the variables $[SSI]$ and $[ISI]$ we are going to use the following approximations or closures \cite{Keeling1999}:
\begin{align*}
[SSI]&\approx\frac{n-1}{n}\cdot \frac{[SS][SI]}{[S]}=\frac{n-1}{n}\cdot \frac{[SS][SI]}{N-[I]}, \\
[ISI]&\approx\frac{n-1}{n}\cdot \frac{[SI]^2}{[S]}=\frac{n-1}{n}\cdot \frac{[SI]^2}{N-[I]},
\end{align*}
where $n(t)$ is the current mean degree of the network,
\[
n(t)=\frac{2[SI]+[SS]+[II]}{N}.
\]
Substituting these into the set of differential equations above we obtain the following approximation:
\begin{subequations}
\begin{align}
\dot{[I]}&=\tau [SI] - \gamma [I],  \label{mastereqa0} \\
\dot{[SI]}&=\gamma ( [II] - [SI])+ \tau \left( \frac{n-1}{n} \right ) \frac{[SS][SI]}{N-[I]} - \\ \nonumber
 &-\tau \left( \frac{n-1}{n} \right) \frac{[SI]^2}{N-[I]}-(\tau +u_1)[SI], \\
\dot{[II]}&=-2\gamma [II] +2\tau \left (\left( \frac{n-1}{n} \right ) \frac{[SI]^2}{N-[I]}+[SI]\right),\\ \label{mastereqd0}
\dot{[SS]}&=2\gamma [SI] - 2\tau \left( \frac{n-1}{n} \right ) \frac{[SS][SI]}{N-[I]} \\ \nonumber
&+ u_2((N-[I])(N-[I]-1)-[SS]).
\end{align}
\end{subequations}

\subsection{Dynamical behaviour}
In the Appendix, we show that the system has two steady states, (1) the disease-free steady state
and (2) the endemic steady state, and the system also exhibits a stable limit cycle as it
undergoes a Hopf bifurcation, see top panel of Fig.~\ref{fig:fig1}.  In the Appendix we also give
the detailed calculations corresponding to the stability analysis. The system is characterised by
three main behaviours as illustrated in Fig.  \ref{fig:fig1} on the $(u_1,u_2)$ parameter plane,
see bottom panel. The first case, from left to right, is when the endemic steady state is stable.
In this case after a short period of damped oscillations, the system settles to the endemic steady
state, for which $[I](t_{final}) \neq 0$. The second case is when both the endemic and the
disease-free steady states are unstable. In this case, the system variables exhibit stable
oscillations, which fail to damp due to the instability of both steady states. Finally, the third
case is when the disease-free steady state is stable, the infection eventually disappears from the
system, and due to the accumulation of the $SS$ edges the network will become fully connected.
Hence, the final state of the system is a complete network with every node being in state $S$. The
curve of transcritical bifurcation is given by Eq.~\eqref{transcrit} (e.g. $u_1=\tau
(N-2)-\gamma$) and the Hopf bifurcation set is defined in Eq.~ \eqref{hopfset}. Obviously, varying
parameters such as $\tau$ will not alter the qualitative behaviour, but the stability regions of
the various steady states change.

A key ingredient to consider in such models is the relation between the dynamics of the epidemic
and the network. The current system can be considered as and adaptive of dynamic network model
\cite{szab2012detailed}, where the epidemic affects link deletion and creation, since these are
type-dependent, and at the same time, link activation deletion can favour or hinder epidemic
spread, respectively. The impact of this interaction is maximal if both processes operate on a
comparable time scale. When this is not the case, the system can exhibit a seemingly surprising
behaviour. In the case of small values for $\tau$ and $u_2$ and a comparably large value for
$u_1$, such that the disease-free steady state is still unstable  (i.e.  $u_1 = \tau(N-2)
-\gamma$), a seemingly eradicated epidemic re-appears at a significant level, see
Fig.~\ref{fig:fig2}. This can be explained as follows.  The low rate of infection combined with a
low rate of link creation, but a high rate of $SI$ edge cutting pushes the system close to the
disease-free steady state, infection slowly disappears from the system. But when the system gets
close to the disease-free state, the cutting of the $SI$ edges is less significant as there are
few such edges and in the mean time the number of $SS$ edges is slowly building up. Hence, a
network that becomes better connected with a very small seed of infection can spark an epidemic
outbreak. Obviously, in a stochastic model it may not be feasible for the system to visit states
of very low prevalence without the epidemic becoming extinct. If we would like to control a system
in this way, it could be quite effective (infection could almost be removed from the system) but
it is crucial to stop or alter the control at the right time, before another epidemic can start.

Concluding our analysis of constant control, we note that there is a wide range of parameter value
combinations that lead to the eradication of the disease. Usually, this requires the deletion of
$SI $ edges at a fast rate at the expense of a dramatic drop in the mean degree of the network.
The system then compensates by connecting susceptible individuals, and in the successful control
case the network becomes completely connected, which is also a dramatic change. Trivially, we can
delete $SI$ edges at a very fast rate, then wait for the infecteds to recover without creating
extra $SS$ links followed by the creation of $SS$ edges in order to reach the desired target
connectivity in the network. Obviously, this strategy requires extremely high cutting rates that
are not feasible in practice, and it will not work in the constant control case and within the
control horizon $T$.  To achieve this type or similar control, in the next section we consider
dynamic control using the Nonlinear Model Predictive Control (NMPC) algorithm.

\section{Time dependent control}
We have seen in the previous section that constant control is not an effective way to control the mean degree of the network and it is a very costly way to control the infection itself.
In particular, cutting infection by breaking down the network is an extreme measure which in reality would correspond to a major quarantine at population level. This
is obviously not feasible, and while the cutting of some potential risky links in response to an epidemic is possible, in general individuals will aim to maintain some form of social connectedness.
Hence, a realistic control should be able to eradicate the disease without leading to a heavily fragmented population. So in this section
we introduce a more sophisticated form of control, i.e. time dependent control.

The basic idea of time dependent control is that we can update the control signal from time to
time according to the current state of the system and our goals. So in this case the control
signals $u_1,u_2$ will be piecewise constant functions. These functions should be bounded by some
realistic values. We want $u_1$ to be positive, since creation of links between infected  and
susceptible individuals would hinder control. But this time we want to admit negative values for
$u_2$ since deleting $SS$ edges will prove useful in controlling the mean degree. There should
exist constants $M_1, M_2$ such that $u_1 \leq M_1$ and $|u_2| \leq M_2$. We introduce a step size
$\Delta t$ for how often we can intervene and change the amount of control, and a constant $T$
which will mark the total length of the control period. We will set a target value for the two
variables we wish to control: $[I^*]$ for the number of infected individuals and $n^*$ for the
mean degree. Using these notations, we can define what we mean by controllability.

\begin{defi}
The system is $\varepsilon$-controllable in time $T$ with step size $\Delta t$ and with control bounds $M_1$, $M_2$ to the targets $[I^*]$, $n^*$, if there are piecewise constant functions
$u_1,u_2:[0,T]\to \mathbb{R}$, such that
\begin{itemize}
\item $0\leq u_1(t) \leq M_1$, $|u_2(t)|\leq M_2$ for all $t\in [0,T]$,
\item $u_1$ and $u_2$ are constants in the intervals $[(k-1)\Delta t,k\Delta t)$ for all $k=1,2,\ldots, [T/\Delta t]$,
\item $|[I](T)-[I^*]|\leq \varepsilon$ and $|n(T) - n^*| \leq \varepsilon$.
\end{itemize}
\end{defi}

Total controllability ($\varepsilon = 0$) would be ideal, but in most cases this is simply too
much to expect from such a control scheme. In practice usually different forms of asymptotic
controllability (termed asymptotic stability) is expected of NMPC algorithms, see \cite{Grune} for
a wide variety of examples. We study finite-time controllability, so our definition using the
error term $ \varepsilon$ is in good agreement with the notion of asymptotic controllability.

We can group the parameters in the following way: the \emph{system parameters} are $N, \tau, \gamma$, $[I](0)$, the \emph{control parameters} are $T, \Delta t, M_1, M_2$, the \emph{targets}
are $[I^*],n^*$ and the error term $\varepsilon$, see Table \ref{paramtable}. Our aim is to investigate how the controllability of the system depends on these parameters.
\begin{table}
\begin{tabular}{ |l|l|l| }
\hline
\multicolumn{3}{ |l| }{\textbf{System parameters}} \\
\hline
$N$ & size of population & $1000$ \\ \hline
$\tau$ & rate of infection across a contact &  \\ \hline
$\gamma$ & rate of recovery & $1$ \\ \hline
$D$ & average length of the infectious period & $1/\gamma$ \\ \hline
$[I](0)$ & number of infecteds at $t=0$  & $10$ \\ \hline
$[SS](0), [SI](0), [II](0)$ & link types at $t=0$ &  \\
\hline
\multicolumn{3}{ |l| }{\textbf{Control parameters}} \\
\hline $u_1$ & rate of cutting $SI$ links & $0 \le u_1 \le M_1$ \\ \hline $u_2$ & rate of
creating/cutting $SS$ links  & $|u_2| \le M_2$ \\ \hline $T$ & time to end of control & $DQ$ \\
\hline $U$ & number of intervention during $D$ &  \\ \hline $\Delta t$ & step size for control
adjustment & $D/U$ \\ \hline
$\varepsilon$ & error term &  \\
\hline
\multicolumn{3}{ |l| }{\textbf{Targets to achieve}} \\
\hline
$[I^{*}]$ & number of infecteds at $T$ & $0$ \\ \hline
$n^{*}$ & average connectivity at $T$ & $n(0)$ \\ \hline
\multicolumn{3}{ |l| }{\textbf{Damping parameters}} \\
\hline
$\lambda_1$ & controlling level of infection &  \\ \hline
$\lambda_2$ & controlling jumps in $u_1$ &  \\ \hline
$\lambda_3$ & controlling average connectivity & \\ \hline
$\lambda_4$ & controlling jumps in $u_2$ & \\ \hline
\end{tabular}
\caption{Table summarising system and control parameters, as well the ideal outcome or target of control. With the applicability in mind, we
work with the average infectious period $D$. Based on this the time to the end of control is set as $T=QD$, and the number of interventions are also
per average length of infection, $\Delta t=D/U$. The control bounds, $M_1$, $M_2$ and $Q$ can also be interpreted as parameters, and will be treated as such.}
\label{paramtable}
\end{table}

We will fix some of the parameters, such as $N=1000$ and $[I](0)=0.01N=10$. Let $D$ be the length of the epidemic, so the recovery rate is $\gamma = \frac{1}{D}$.
Let $Q > 0$ be a constant such that $T=D \cdot Q$, meaning that we can set control over many generations/waves of infection. We also make the frequency of intervention of control to depend on
$D$, and set $U$ to be the parameter for how many times we are to intervene during an average infectious period, so for the step size for control we use $\Delta t=\frac{D}{U}$.
For simulation purposes we used $D=1$, $Q=10$ and $U=5$. This means a control period of $T=10$ and a step size of $\Delta t=0.2$. While these are in arbitrary units, these values translate to
seeing an intervention every day or every week for disease with a typical average infectious period of 5 days or 5 weeks, respectively. We also have to provide some reasonable values for $M_1$
and $M_2$. For example, if
\[
u_1 \cdot \Delta t = 0.2,
\]
then this corresponds to deleting $20\%$ of the $SI$ edges in $\Delta t$ time. This is quite a considerable amount, and hence,  the maximum value $M_1$ for $u_1$ is set to $\frac{0.2}{\Delta t}=
\frac{0.2 \cdot U}{D}$, which for our simulation parameters equates to $1$. Similarly, an appropriate value for $u_2$ is $\frac{u_1}{N}=0.001$, since $u_2$ has a quadratic multiplier in terms of $N
$ in the system of equations \eqref{mastereqa0}-\eqref{mastereqd0},
while the multiplier of $u_1$ is linear in $N$.

Our targets will be $[I^*]=0$ and $n^*=n(0)$ describing our goal that we wish to find and apply a control which eradicates infection while keeping the network connected. In this case without
loss of generality we set the target average connectivity to its value at time $t=0$. Finally, the error will be acceptable if it is lower than $0.1$, but ideally it should be of much smaller
magnitude than this value.  Nonetheless, we say the control is effective if $\varepsilon \leq 0.1$.

\subsection{Nonlinear model predictive control}
Here, for the readers convenience, a brief introduction to Nonlinear Model Predictive Control
(NMPC) is provided. NMPC is a control strategy which is suited for constrained, multivariable
problems. The main idea of the method is as follows. At each step of the NMPC algorithm a sequence
of optimal control signals is calculated along a prediction horizon of fixed length by minimizing
an objective functional which includes predicted future outputs of the system. This optimization
is a nonlinear programming problem which is solved subject to some constraints imposed on the
input and output signals. Only the first control of the obtained sequence of optimal signals is
applied to the system, then the prediction horizon is moved one step forward and the next control
signal is calculated the same way. Due to this moving horizon technique the NMPC is also called
Receding Horizon Control. There are many applications of NMPC, for example controlling drug
dosing, industrial plants or automobiles, see the collection of survey papers \cite{Magni}. For
further theoretical details on NMPC, we refer to the monograph \cite{Grune}.

Our aim is now to apply the NMPC method to control epidemics spread. We use a little different system than in the previous sections:
\begin{align} \label{mastereqn0}
\dot{[I]}&=\tau [SI] - \gamma [I], \\
\dot{[SI]}&=\gamma ( [II] - [SI])+ \tau ([SSI]-[ISI]-[SI])-u_1[SI], \\
\dot{[II]}&=-2\gamma [II] +2\tau ([ISI]+[SI]), \\  \label{mastereqn1}
\begin{split}
\dot{[SS]}&=\gamma [SI] - 2\tau [SSI]+\max\{u_2,0\}\cdot ((N-[I])(N-[I]-1)-[SS]) \\
&\quad{}+\min\{u_2,0\}\cdot [SS].
\end{split}
\end{align}
We now admit the algorithm to assign negative values to $u_2$, so that it can also delete $SS$
edges, equation \eqref{mastereqn1} is adjusted accordingly. The vector of state variables and
control variables will be $x=([I], [IS], [II], [SS])$ and $u=(u_1, u_2)$, respectively. The output
variables are the number of infected individuals and the mean degree, so $y=([I], n)$. The $i$th
coordinate of $x$ and $y$ will be denoted by $x_i$, $y_i$, respectively, e.g., $y_1=[I]$ and
$y_2=n$.

In order to apply the NMPC algorithm, we should first discretize the system
\eqref{mastereqn0}--\eqref{mastereqn1}. We fix a time step $\Delta t$ and observe the system only
at instants $t=k\Delta t$ where $k\in\mathbb{Z}$. For simplicity, we shall omit $\Delta t$ and
write $x(k)$, $y(k)$ which means that $x$ and $y$ are evaluated at time instant $k\Delta t$. We
suppose that the control variables $u_1$ and $u_2$ are held constant along the intervals $[k\Delta
t,(k+1)\Delta t)$ ($k\in\mathbb{Z}$), in other words they are piecewise constant functions. With
these conventions we obtain the following discretized form of system
\eqref{mastereqn0}--\eqref{mastereqn1}:
\begin{align} \label{nmpc0}
x(k+1)&=F(x(k),u(k)), \\
y(k+1)&=h(x(k+1)), \label{nmpc1}
\end{align}
where $x(k)\in \mathbb{R}^4$ is the vector of state variables, $u\in\mathbb{R}^2$ is the vector of
input (control) signals and $y\in \mathbb{R}^2$ is the vector of output signals. Furthermore, the
function $F$ symbolizes that we solve the system of ODEs given by
Eqs.~\eqref{mastereqn0}--\eqref{mastereqn1} numerically on the interval $[k\Delta t,(k+1)\Delta
t]$ and $h(x_1,x_2,x_3,x_4)=(x_1,(2x_2+x_3+x_4)/N)$. As we explained before, we impose the
following constraints on the control signals:
\begin{align*}
0 &\leq u_1(k) \leq M_1, \\
-M_2 &\leq u_2(k) \leq M_2.
\end{align*}
Now, the control action $u$ at time $k$ is computed as follows. We fix a prediction horizon of
length $P$ steps and perform a nonlinear optimization procedure over the admissible set of future
control actions as described below. Denote by $u_i(k+j|k)$ ($i=1,2$) an arbitrary admissible
future control action at time $k+j$ chosen at time $k$. If we choose admissible sequences of
future control actions $u_i(k|k), u_i(k+1|k),\dots, u_i(k+P-1|k)$ ($i=1,2$), then these controls
yield predicted future outputs $y_1(k+j|k)$, $y_2(k+j|k)$ ($j=1,2,\dots,P-1$), where the notation
means that $y_i(k+j|k)$ is a predicted output at time $k+j$ calculated at instant $k$. More
specifically,
\begin{align*}
x(k+j|k)&=F\bigl(x(k+j-1|k),(u_1(k+j-1|k),u_2(k+j-1|k))\bigr),\\
y(k+j|k)&=h(x(k+j|k)),\qquad j=0,1,2,\dots,P-1,
\end{align*}
where $x(k+j|k)$ denotes the predicted state at $k+j$ calculated at instant $k$. The setpoints
for the output signals are $y_{1_s}=[I^*]=0$ and $y_{2_s}=n^*=n(0)$, therefore, we choose the
objective functional $J\colon\mathbb{R}^{2P}\to\mathbb{R}$ to be minimized to have the form
\begin{align*}
J\bigl(u(k|k),\dots,u(k+P-1|k)\bigr)&=\sum_{j=0}^{P-1} \lambda_1 (y_1(k+j|k))^2+\lambda_2 (\Delta u_1(k+j|k))^2\\&\quad{}+\lambda_3 (y_2(k+j|k)-n(0))^2+\lambda_4 (\Delta u_2(k+j|k))^2
\end{align*}
with parameters $\lambda_1,\dots,\lambda_4$ where $\Delta u_i(k+j|k)=u_i(k+j|k)-u_i(k+j-1|k)$ is
the  predicted control effort at instant $k+j$ calculated at time $k$. The evaluation of the
functional $J$ requires the numerical solution of the ODE system give by
Eqs.~\eqref{mastereqn0}-\eqref{mastereqn1}. Clearly, by adjusting the parameters we can put more
weight on the quadratic difference terms or penalize large control efforts that will be discussed
later in details. For example, $\lambda_1$ penalises small departures from the no-epidemic state,
while $\lambda_2$ penalises large changes in rewiring rates. Now performing the above nonlinear
optimization problem we obtain a sequence of optimal controls $u_i(k|k), u_i(k+1|k),\dots,
u_i(k+P-1|k)$ ($i=1,2$). This can be done by using a nonlinear optimisation routine, such as
\textit{lsqnonlin} in \textit{Matlab}, with the quadratic functional $J$ as an input. Then only
$u(k):=(u_1(k|k),u_2(k|k))$ is applied to the system and the prediction horizon is translated one
step forward and the same optimization procedure is implemented to calculate the next control.

It is intuitively clear that dynamical control is more effective compared to constant control. However, the number of parameters involved in setting up or specifying
dynamic control makes it non-trivial to understand which combinations of factors and what parameter values will make the system controllable. In the next section, we will
numerically explore in detail the impact of system, control and damping parameters, see Table~\ref{paramtable}.

\subsection{The interplay between the infection rate and control bounds}

First, let us analyze the prevalence level, i.e. the number of infected individuals at time
$T=10$, $[I](T)$, and the  mean degree at time $T$ in the uncontrolled system ($u_1, u_2 \equiv
0$) for different values of $\tau$. As expected, in Fig. \ref{fig:fig3} we can see that for very
small values of $\tau$ (for about $\tau < 0.05$) the infection disappears from the system even
without control. However, for higher values of $\tau$ the disease becomes more widespread and the
prevalence level converges towards the full population size. When no control is applied the mean
degree of the system remains unchanged in each step, so naturally the final value of the mean
degree is the initial value $n(0)=10$ for each $\tau$.

Now using the NMPC introduced above, the case of dynamic control is studied. Initially, we
consider a set of fixed control parameters, $M_1$ and $M_2$,  and a varying value of $\tau$.
Naturally, it is easier to control the infection when the infection rate is low and impossible
within the given control bounds if the infection rate is high. In Fig.~\ref{fig:fig3}, the
prevalence and the mean degree at time $T$ for different values of $\tau$ are plotted. The figure
shows that for approximatively $\tau > 0.15$ the control is ineffective, both the value of
$[I](T)$ and $n(T)$ visibly differ from their target at the end of the control period. For about
$\tau >0.25 $, the final number of $[I]$ is greater than $[I](0)=10$, so in this case the control
failed to decrease the initial amount of infected individuals. For even higher values of $\tau$,
the control has little effect on $[I]$ or $n$, so if $\tau \to \infty$ the final values of these
variables converge to the final values of an uncontrolled system.  The final mean degree of 10 can
be attained for some small vales of $\tau$ which cannot be said in the uncontrolled case. However,
the value of $n(10)$ becomes much lower for higher vales of $\tau$ despite the control.

This behaviour is due to the strict bounds on the values of the control parameters: limiting the
cutting rate of $SI$ edges, $u_1$, throughout the entire control period makes the control
inefficient for higher values of $\tau$. For high infection rates even $SS$ edges are cut, but
again with a limited strength and making little difference. In fact, this only results in the drop
of the mean degree, since in this case $u_2$ has no capacity to make new connections. To shed some
light on the precise dependency of successful control on the bounds of the rewiring rates for
different values of $\tau$, a detailed numerical exploration is carried out.  To carry out this
exploration, we fix the damping parameters as follows $\lambda_1 = 10^4$ and
$\lambda_2=\lambda_3=\lambda_4=1$. This choice penalises even a small departure from the ultimate
target of disease eradication. While, we fix these, the damping parameters themselves will impact
on the controllability of the system, and this is considered in the next subsection.

First, we will investigate the effect of $M_1$'s magnitude with fixed values for $\tau$ and $M_2$.
A value of $\tau=2$ is a good starting point given that with the previous bounds for link
rewiring, $M_1=1$ and $M_2=0.001$, control was not successful even for $\tau=1$, see
Fig.~\ref{fig:fig3}. If we wish to keep $M_2=0.001$, we should increase the value of $M_1$. Figure
\ref{fig:fig4} shows that $M_1=18$ makes the system controllable.

Extensive numerical simulations suggest that for a fixed value of $\tau$ and $M_2$, there is a
critical value $M_1^c$, such that if $M_1$ is lower than $M_1^c$, the control is not effective.
However, if  $M_1$ is larger than the critical value, then control is effective in $T$ units of
time. The higher the value of $M_1$ the less time is needed to control the system. But choosing a
high value for $M_1$ implies that control is more severe or drastic. Hence, if our aim is to
control our system in $T$ units of time by using the least invasive control, it is optimal to
choose $M_1^c$ as the bound for $M_1$. This critical value is the strictest bound admissible. In
Fig.~\ref{fig:fig5} (left panel) the critical value $M_1^c$ for three different values of $M_2$ is
plotted as $\tau$ is varied. These curves in fact define the strictest possible bounds, and hence,
one can use these to identify $(M_1,M_2)$ pairs that can deliver a successful control. Moreover,
the same figures shows that higher values of $M_2$ have negligible effect on the critical $M_1$
curve, since the fast creation of $SS$ does not help to control the epidemic. In Figure
\ref{fig:fig5} (right panel) the critical value of $M_1$ is plotted for a range of $M_2$ values
and different infection rates. The same applies as previously: choosing bounds below this curve
will not result in an effective control. Choosing a pair $(M_1,M_2)$ belonging to these curves is
in some sense optimal, since these represent the strictest bounds.

\subsection{Effects of $\Delta t$ and the damping parameters on controllability}

In this section we analyze how the value of the step size $\Delta t$ and the damping parameters
(i.e. $\lambda_i$ $(i=1,2,3,4)$) in the cost functional affect system controllability. Let us
first deal with the step size. A greater value for this means a slower reaction, so as we increase
it, controlling the system requires more radical changes in control, and the change in the mean
degree during the control period could be quite drastic. However, we experienced that step sizes
$\Delta t \leq 5$ are effective - which means $U=0.2$ (i.e. $\Delta t=D/U=1/0.2=5$) is not enough,
but any larger $U$ suffices (the parameter $U$ marked the number of control actions during the
average infectious period $D$, and $U$ needs to take values less than one if one wants to
investigate slow reactions in control). For a greater step size, the reaction of the control is
not fast enough to control the system in $T=10$ units of time. Figure \ref{fig:fig7} uses
$M_2=0.5$, $\tau=1$ and the critical value of $M_1^c=7.8$. In Fig.~\ref{fig:fig7} the effect of
control is shown for four different values of $\Delta t$. It is clear that the system is only
controllable if time steps are small enough. While we do not separately investigate the effect of
the control parameter $T$, we note that an increase in the control horizon is likely to make
controllability possible.

Let us continue with the analysis of the damping parameters. The damping parameters assigned to
$\Delta u_1$ and $\Delta u_2$ are $\lambda_2$ and $\lambda_4$, respectively. When both are large
compared to $\lambda_1$ and $\lambda_3$, achieving the control target will be difficult due to
small increments in the rewiring rates. As shown in Fig.~\ref{fig:fig8} (continuous line),
infection seems to eradicated after a large excursion into high infection levels, but network
connectivity is far from the target. This is exacerbated by a magnitude difference in size of the
$\lambda_1$ and $\lambda_3$, with control focused more on achieving eradication of the disease.
However, when the control functional depends solely on controlling the spread, then this target is
quickly achieved, but this happens at the price of the network being completely disconnected, see
Fig.~\ref{fig:fig8} (dashed line). When, the control adjustment is not penalised and with a
stronger focus on achieving the target connectivity, the system proves to be uncontrollable since
the disease cannot be eradicated at the end of the control period, see Fig.~\ref{fig:fig8} (dotted
line).

\subsection{Control-bound-induced targets}

Posing a controllability question usually involves establishing the control bounds for a given
target. However, understanding what targets can be achieved with given control bounds is equally
valuable, especially when these could be close to the ideal targets. We have seen in the previous
sections that if we fix a value of the constraint (i.e. $M_2$) on $u_2$ and the infection rate
$\tau$, there exists a critical value for $M_1$ below which the system is not controllable. In
many cases, the main difficulty was posed by reaching the target connectivity. More importantly,
the infection is almost completely eradicated from the system in every case, if the formerly fixed
$\lambda_1=10^4$ and $\lambda_2=\lambda_3=\lambda_4=1$ damping parameters are used. So for a
weakened control, let us admit a decrease in the value of the target mean degree. For example, if
we use the previously seen $M_2=0.5$, $ \tau=1$ parameters, we have seen that the critical value
$M_1^c$ was $7.8$, and Fig.~ \ref{fig:fig5} (left panel) shows that the system is not controllable
for $M_1=6$. Now let us use the target value $n^*=7.5$ admitting a $25\%$ decrease in the mean
degree. In Fig.~\ref{fig:fig6}, it is clearly illustrated that the $n^*=10$ cannot be achieved,
see top row. However, modifying the target to $n^*=7.5$, the system becomes controllable, see
bottom row. Let us fix the parameters $M_2$ and $\tau$ above and analyse the highest possible
achievable $n^*$ for different values of $M_1$. Table \ref{table2} below shows the results of some
of our simulations.

\begin{table}
\begin{center}
  \begin{tabular}{  p{2.5cm}|  p{2.5cm}  }
    $M_1$ & $n^*$ \\ \hline \hline
    7.8 & 10 \\ \hline
    7.5 & 9.2\\ \hline
    7 & 8.6  \\ \hline
    6.5 & 8.2  \\ \hline
    6 & 7.6  \\ \hline
    5.5 & 7.2  \\ \hline
    5 & 6.6 \\ \hline
    4.5 & 6\\  \hline
    4 & 5.2 \\ \hline
    3.5 & 4.4 \\
  \end{tabular}
\end{center}
  \caption{Table showing the achievable target $n^*$ for different values of the control bound $M_1$.}
    \label{table2}
\end{table}

\section{Discussion}

The control in this paper does not appear in the form of what could be termed as classic control.
More precisely, classic control problems in epidemiology involve the minimisation of an integral or cost function.
Here, we focus on the end target and we select the piece-wise constant control signal that allows us to be as close as
possible to the final target. Obviously, within the control parameters that we assume, we ignore costs and a cumulative measure
of the amount of intervention and costs due to infection. This can obviously be built into further models.

While setting up the problem in this way has been a first step to bridge the gap between classic
compartmental control and modern epidemiological models, it is straightforward to apply the same
methodology to more complicated settings involving costs and competing effects, such as the trade
of in cost between vaccination and the number of infectious individuals, namely more vaccination
increase the cost, but results in less infectious cases, which in turn reduces cost. In our case
this trade of was realised by aiming to control disease spread while maintaining social cohesion.
Obviously, if the network cohesion is not required, control will lead to the trivial case of
cutting the network to the point where transmission is no longer possible. In real life this is
not the case, as for STIs persist due to the network being well connected with many concurrent
partnerships, and it is reasonable to assume that control will need to be achieved without
breaking the network of contacts completely.

The next step for this method is to extend it to individual-based network simulations, and work
out to what extent the control predicted by the pairwise model would also translate to
good/optimal control in the stochastic network model. Such endeavours already exist and the first
signs are positive in that control from mean-field type models seem to translate, at least for
some cases, to the simulation counterpart \cite{JC_KAJW_KT_2013}.

\section*{Acknowledgements} P\'{e}ter L. Simon acknowledges support from OTKA (grant no. 81403).

\section*{Appendix: steady states and their stability for the constant control case.}
\label{StabilityAnalysis}

Let us calculate the steady states of system \eqref{mastereqa0}-\eqref{mastereqd0}. These are the solutions of
\begin{subequations}
\begin{align}
0&=\tau [SI] - \gamma [I], \label{mastereqa} \\ \label{mastereqb}
0&=\gamma ( [II] - [SI])+ \tau \left( 1- \frac{N}{2 \cdot[SI] + [II] + [SS]} \right ) \frac{[SS][SI]}{N-[I]} - \\ \nonumber
 &-\tau \left( 1-\frac{N}{2 \cdot[SI] + [II] + [SS]} \right) \frac{[SI]^2}{N-[I]}-(\tau +u_1)[SI], \\
0&=-2\gamma [II] +2\tau \left (\left( 1- \frac{N}{2 \cdot[SI] + [II] + [SS]} \right ) \frac{[SI]^2}{N-[I]}+[SI]\right),
\label{mastereqc} \\ \label{mastereqd}
0&=2\gamma [SI] - 2\tau \left( 1- \frac{N}{2 \cdot[SI] + [II] + [SS]} \right ) \frac{[SS][SI]}{N-[I]} \\ \nonumber
&+ u_2((N-[I])(N-[I]-1)-[SS]).
\end{align}
\end{subequations}
\noindent By solution, we mean an all-real, all-positive solution. It is easy to see that the disease-free steady state of the system is
\begin{align*}
[I]&=0, \\
[SI]&=0, \\
[II]&=0, \\
[SS]&=N(N-1).
\end{align*}
Denoting the disease-free steady state as $E_d$, the Jacobian at state $E_d$  is
\[
J(E_d) =
 \begin{pmatrix}
  -\gamma & \tau & 0 & 0 \\
  0 & -\gamma + \tau (N-2) - (\tau +u_1) & \gamma & 0 \\
  0  & 2\tau  & -2\gamma & 0  \\
  -u_2(2N+1) & 2\gamma-2\tau(N-2) & 0 & -u_2
 \end{pmatrix}.
\]
It is clear that $-\gamma$ and $-u_2$ are eigenvalues of the Jacobian, and these eigenvalues are
always real and negative. So we only have to deal with the eigenvalues of the inner $2 \times 2$
submatrix:
\[
 \begin{pmatrix}
   -\gamma + \tau (N-2) - (\tau +u_1) & \gamma \\
   2\tau  & -2\gamma \\
 \end{pmatrix}.
\]
The determinant of this submatrix is $2\gamma(\gamma-\tau(N-2)+u_1)$, its trace is $-3\gamma+\tau(N-3)-u_1$. For stability we need the eigenvalues to have negative real parts. For this the
determinant has to be positive and the trace has to be negative. So if $u_1 > \tau(N-2)-\gamma$ and $u_1 > \tau(N-3) - 3\gamma$ the disease-free steady state is stable. Note that the second
condition bears no new information, so we can exclude that. Thus our only criterion for the disease-free steady state to be stable is:
\begin{equation}
u_1 > \tau(N-2)-\gamma \label{transcrit}
\end{equation}
Note that in the disease-free steady state, the mean degree is $n=\frac{N(N-1)}{N}=N-1$, so the network becomes fully connected.

To calculate the endemic steady state(s), we first express the variable $[SI]$ from equation \eqref{mastereqa} to get
\[
[SI]=\frac{\gamma}{\tau}[I].
\]
Then we express $[SS]$ from equations \eqref{mastereqb}-\eqref{mastereqd}:
\[
[SS]=u_2(N-[I])(N-[I]-1)-\frac{u_1}{u_2}\cdot \frac{\gamma}{\tau}[I].
\]
We substitute these expressions of $[SI]$ and $[SS]$ in Eq. \eqref{mastereqc}. We obtain a quadratic equation for $[II]$, from which $[II]$ can also be expressed in terms of $[I]$:
\[
[II]=\frac{1}{2}\frac{(A-B\cdot C + \sqrt{D})}{B},
\]
where
\begin{align*}
A&=\frac{[SI]+N-[I]}{N[SI]}, \\
B&=\gamma \cdot \frac{N-[I]}{\tau N[SI]^2}, \\
C&=2[SI]+[SS], \\
D&=(A-B\cdot C)^2-4B(1-A \cdot C).
\end{align*}

Since we are looking for all-real solutions, if $D < 0$, we are left without a solution having a
meaning to us. Otherwise, we substitute these expressions for $[SI]$, $[SS]$ and $[II]$ to
equation \eqref{mastereqb}, and we get an equation containing only the unknown $[I]$. Due to its
complexity, we refrain from writing it out in detail, but let us denote it as equation $(*)$. By
solving equation $(*)$ we get the endemic steady states. Our numerical experiments show, that
there is always only one all-positive solution, so we can conclude that the endemic steady state
(if it exists) is unique. Let us denote this state by $E_e(u_1,u_2)$. The Jacobian is far more
complicated this time, we exclude its concrete form. Substituting $E_e(u_1,u_2)$ into the
Jacobian, we can see by numerical experiments that for some values of $u_1$ there exists a value
$u_2^*$, such that for a lower value of $u_2$ than this $u_2^*$ the Jacobian at $E_e(u_1,u_2)$ has
two real, negative eigenvalues and two imaginary eigenvalues with positive real parts.  For $u_2 >
u_2^* $, the real part of the two imaginary eigenvalues becomes negative. To calculate the exact
value of this $u_2^*$, let us use the method introduced in \cite{szab2012detailed} and write the
characteristic polynomial of the Jacobian at $E_e$ in the following form:
\[
\lambda^4-b_3\lambda^3 +b_2\lambda^2-b_1\lambda +b_0,
\]
such that $b_3 = \Tr J(E_e)$, $b_0 = \det J(E_e)$ and $b_1,b_2$ can be given as the sum of some subdeterminants
of the Jacobian, the concrete form of which is not important at this moment. In the case of $4\times 4$ matrices the necessary and sufficient condition for the existence of pure imaginary
eigenvalues is
\begin{equation} \label{hopf}
b_0b^2_3= b_1(b_2b_3-b_1)\text{ and }\sign b_1 = \sign b_3,
\end{equation}
Thus the Hopf-bifurcation set can be defined as
\begin{equation}
H = \{(u_1,u_2) \in \mathbb{R}^2_+ : \exists \text{ } [I] \in [0,N] \text{ such that $(*)$, \eqref{hopf} hold} \} \label{hopfset}
\end{equation}
This is a simple curve in the $(u_1,u_2)$-parameter plane. $E_e$ is stable above the curve and is unstable below. There is notable oscillation in the value of $[I]$ according to time in the unstable
region.

\newpage

\newpage
\begin{figure}[h!]
        \centering
        \begin{subfigure}[b]{0.8\textwidth}
                \centering
            \includegraphics[width=\textwidth, height=0.24\textheight]{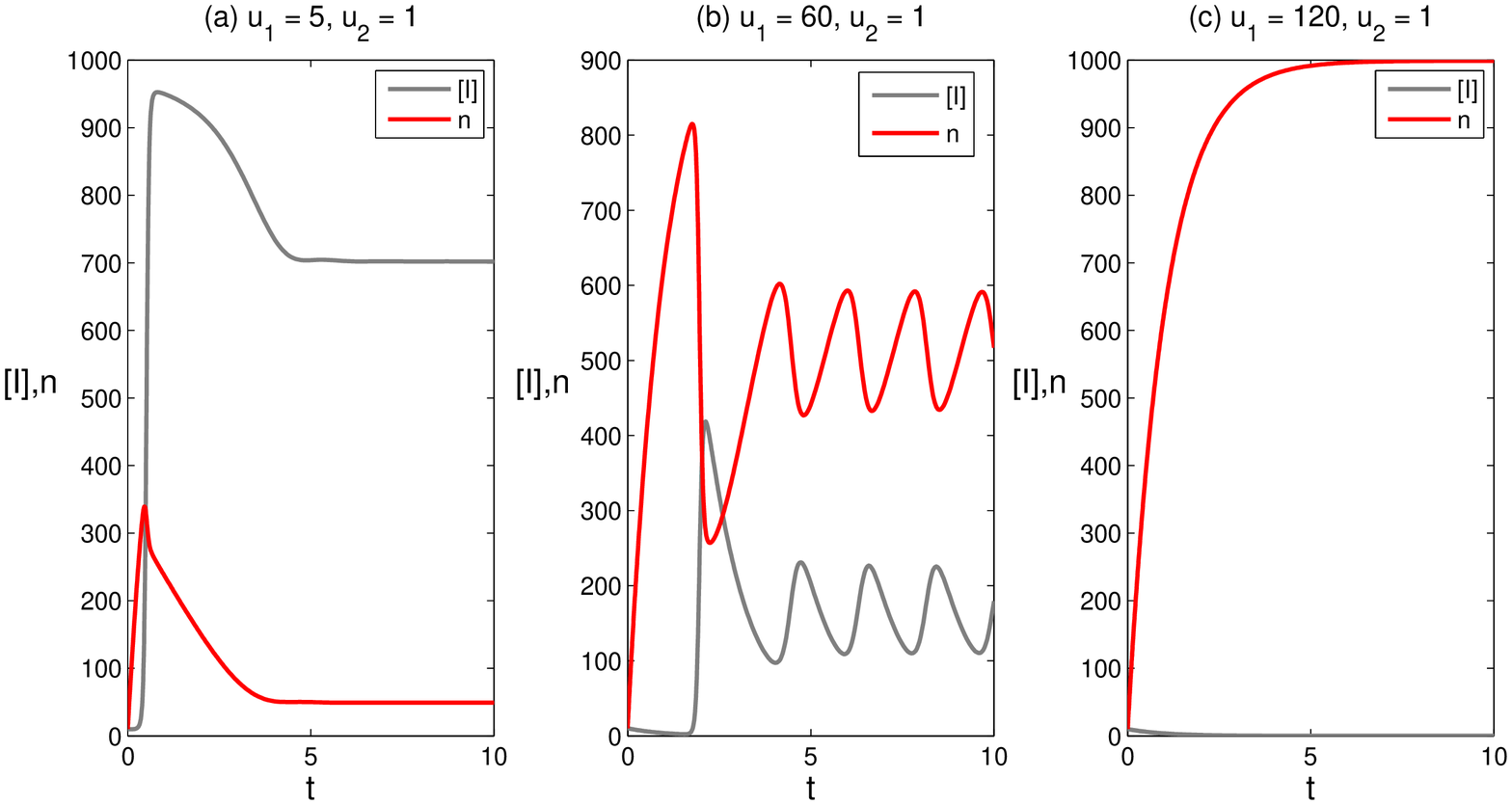}
        \end{subfigure}
        \begin{subfigure}[b]{0.45\textwidth}
                \centering
                        \includegraphics[width=\textwidth,height=0.2\textheight]{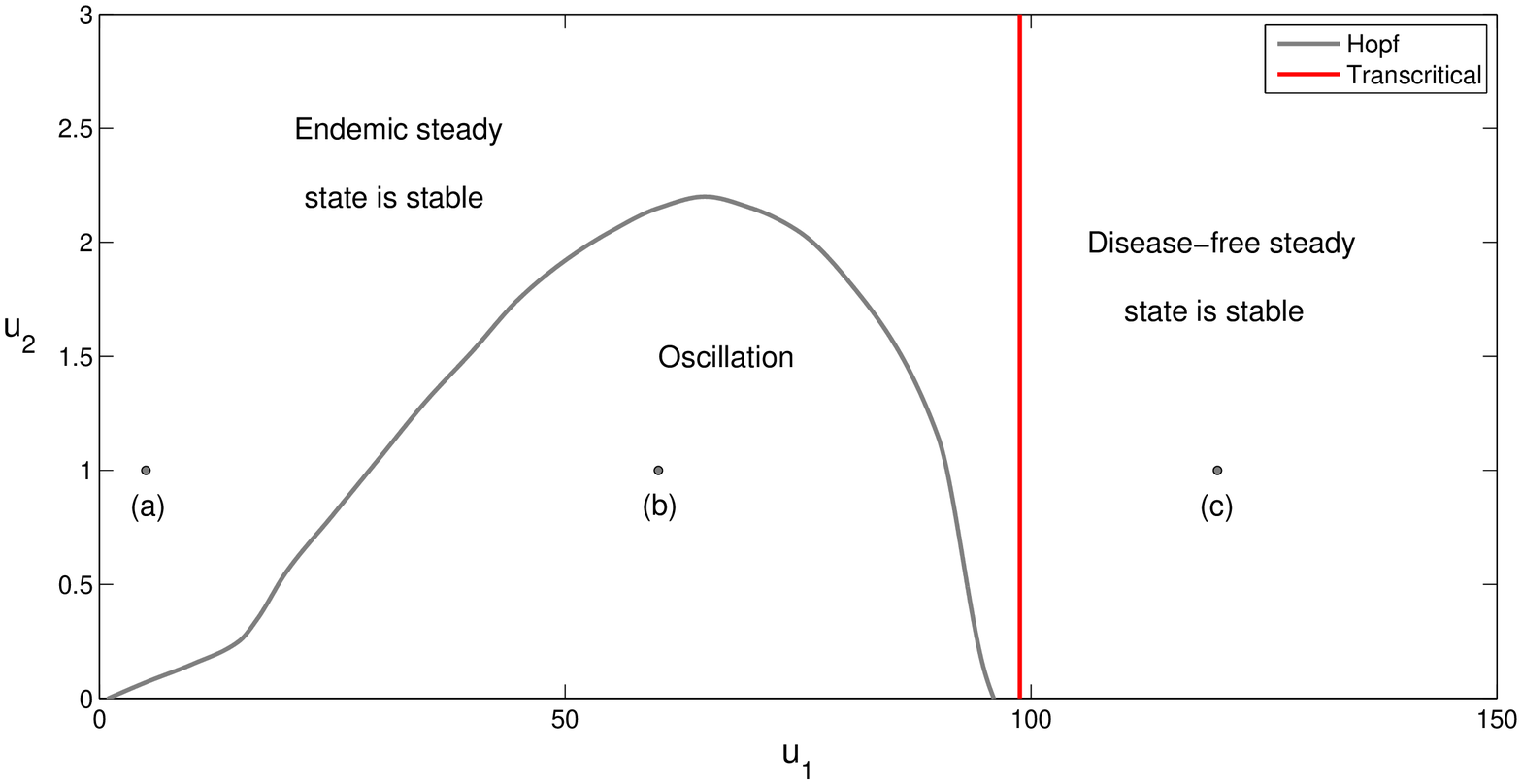}
        \end{subfigure}
        \caption{Typical system behaviours (top row) and bifurcation diagram (bottom panel) for $N=1000$, $n(0)=10$, $\tau=0.1$, $\gamma=1$ and $I(0)=10$.}\label{fig:fig1}
\end{figure}

\begin{figure}[h!]
  \centering
    \includegraphics[width=0.6\textwidth, height=0.3\textheight]{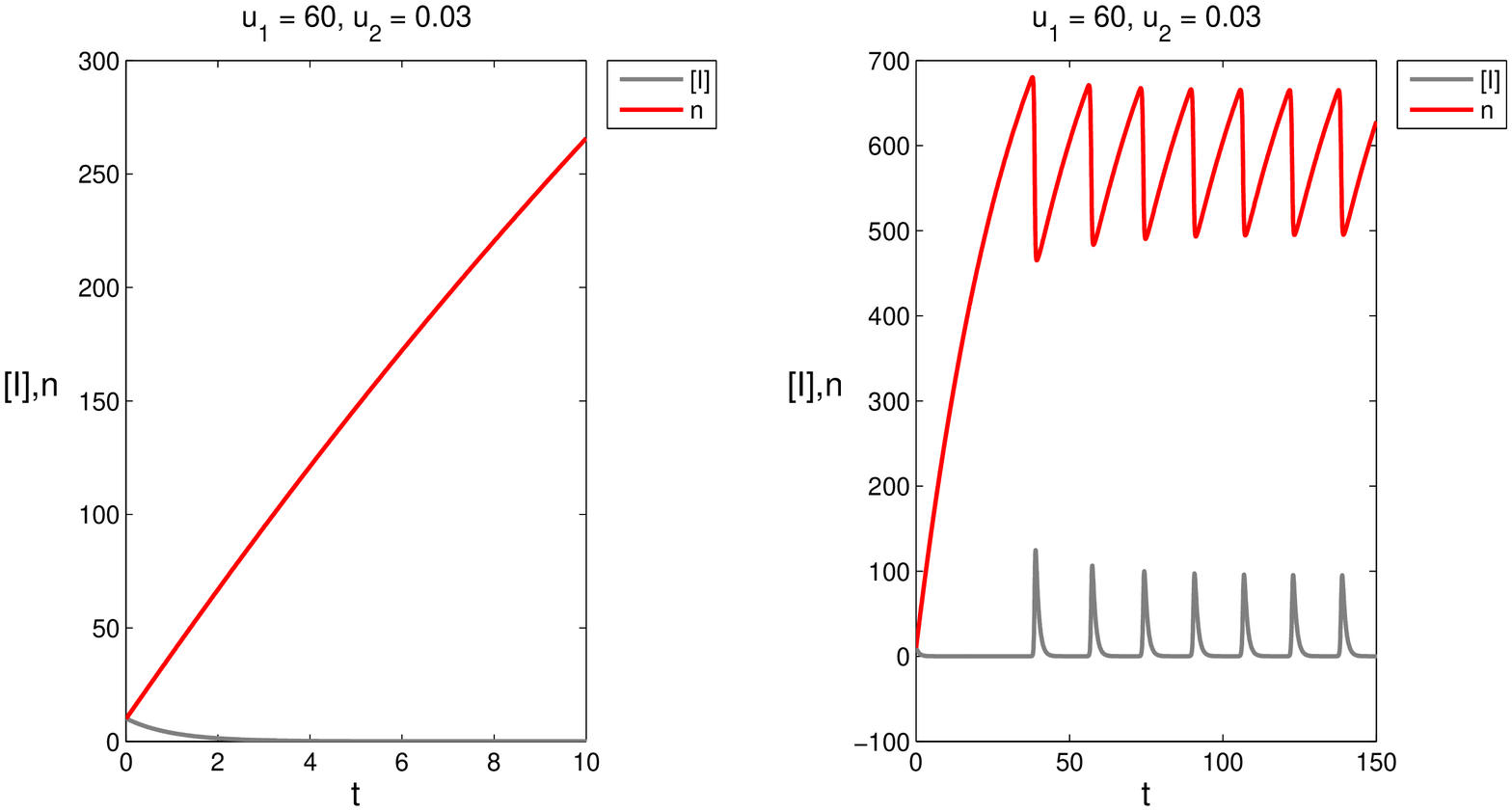}
    \caption{Time evolution of prevalence and network connectivity for $N=1000$, $n(0)=10$, $\tau=0.1$, $\gamma=1$ and $I(0)=10$.} \label{fig:fig2}
\end{figure}

\begin{figure}[h]
  \centering
    \includegraphics[width=0.6\textwidth, height=0.3\textheight]{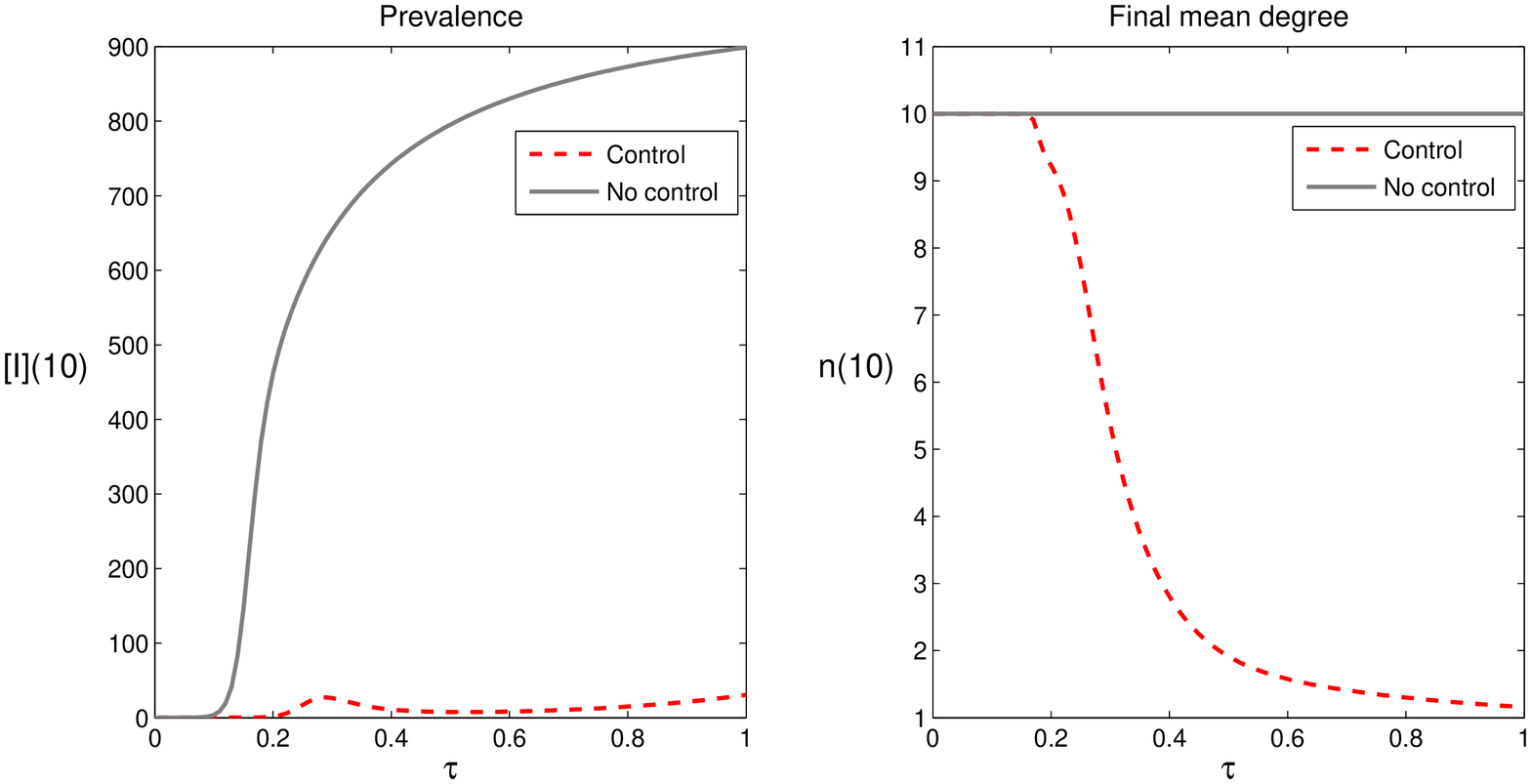}
    \caption{The value of prevalence and network connectivity at the end of the control at time $T=10$, $[I](T)$ and $n(T)$, as a function of the transmission parameter $\tau$ for $\gamma=1$, $N=1000$, $I(0)=10$, $I^*=0$, $n(0)=n^*=10$, $M_1=1$, $M_2=0.001$, $\Delta t=0.1$, $\lambda_1=10^4$ and $\lambda_2=\lambda_3=\lambda_4=1$.} \label{fig:fig3}
\end{figure}

\begin{figure}[h]
  \centering
    \includegraphics[width=0.6\textwidth, height=0.3\textheight]{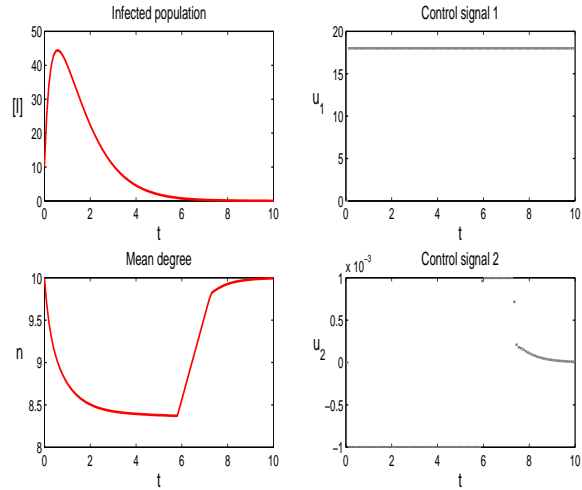}
    \caption{Time evolution of prevalence, network connectivity and control signals, $u_1$ and $u_2$, for $\tau=2$, $\gamma=1$, $N=1000$, $I(0)=10$, $I^*=0$, $n(0)=n^*=10$, $M_1=18$, $M_2=0.001$, $\Delta t=0.1$, $\lambda_1=10^4$ and $\lambda_2=\lambda_3=\lambda_4=1$.} \label{fig:fig4}
\end{figure}

\begin{figure}[h]
  \centering
    \includegraphics[width=\textwidth, height=0.3\textheight]{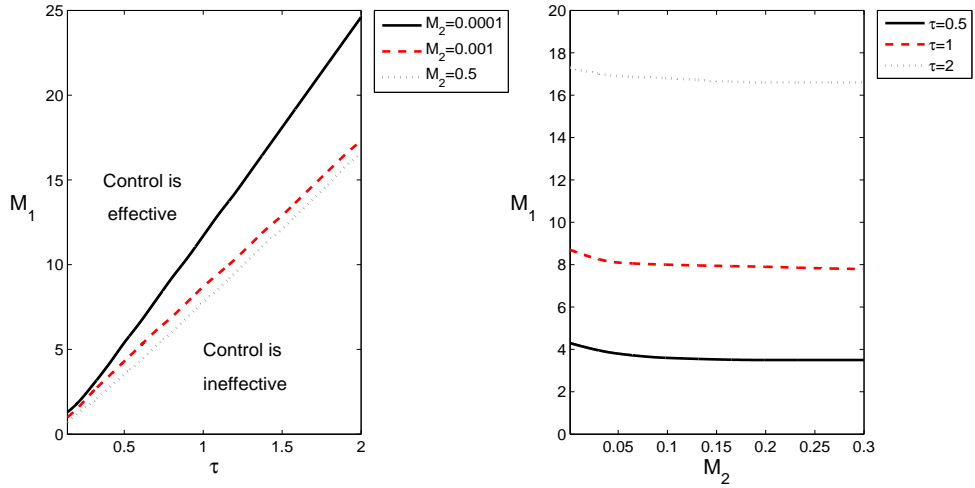}
    \caption{Threshold plots illustrating the relation between system ($\tau$) and control parameters ($M_1$ and $M_2$) for $\gamma=1$, $N=1000$, $I(0)=10$, $I^*=0$, $n(0)=n^*=10$, $\Delta t=0.1$, $\lambda_1=10^4$ and $\lambda_2=\lambda_3=\lambda_4=1$.} \label{fig:fig5}
\end{figure}

\begin{figure}[h]
  \centering
    \includegraphics[width=\textwidth, height=0.3\textheight]{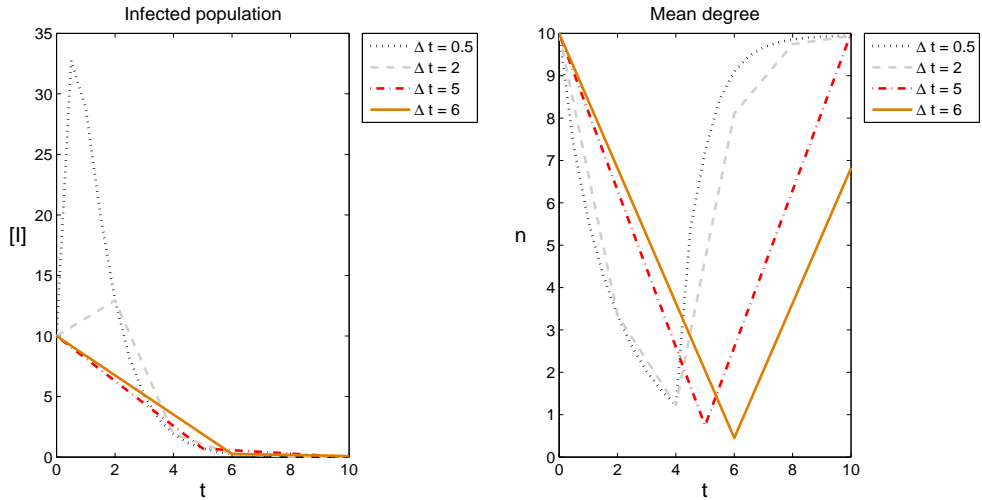}
    \caption{The impact of intervention frequency in terms of the time evolution of prevalence  and network connectivity for $\tau=1$, $\gamma=1$, $N=1000$, $I(0)=10$, $I^*=0$, $n(0)=n^*=10$, $M_1=7.8$, $M_2=0.5$, $\lambda_1=10^4$ and $\lambda_2=\lambda_3=\lambda_4=1$.} \label{fig:fig7}
\end{figure}

\begin{figure}[h]
  \centering
    \includegraphics[width=0.8\textwidth, height=0.3\textheight]{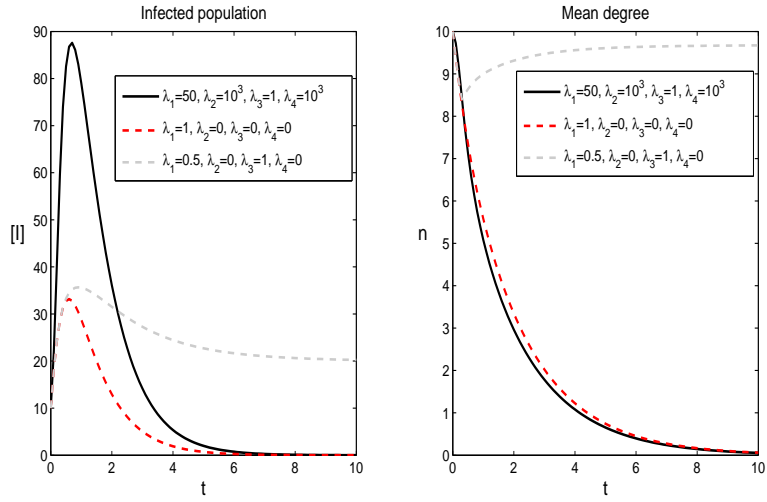}
    \caption{Dependence of system's controllability on the damping parameters in terms of the time evolution of prevalence  and network connectivity for $\tau=1$, $\gamma=1$, $N=1000$, $I(0)=10$, $I^*=0$, $n(0)=n^*=10$, $M_1=7.8$, $M_2=0.5$, $\Delta t=0.1$ and three sets of damping parameters.} \label{fig:fig8}
\end{figure}

\begin{figure}[h]
  \centering
    \includegraphics[width=0.6\textwidth, height=0.3\textheight]{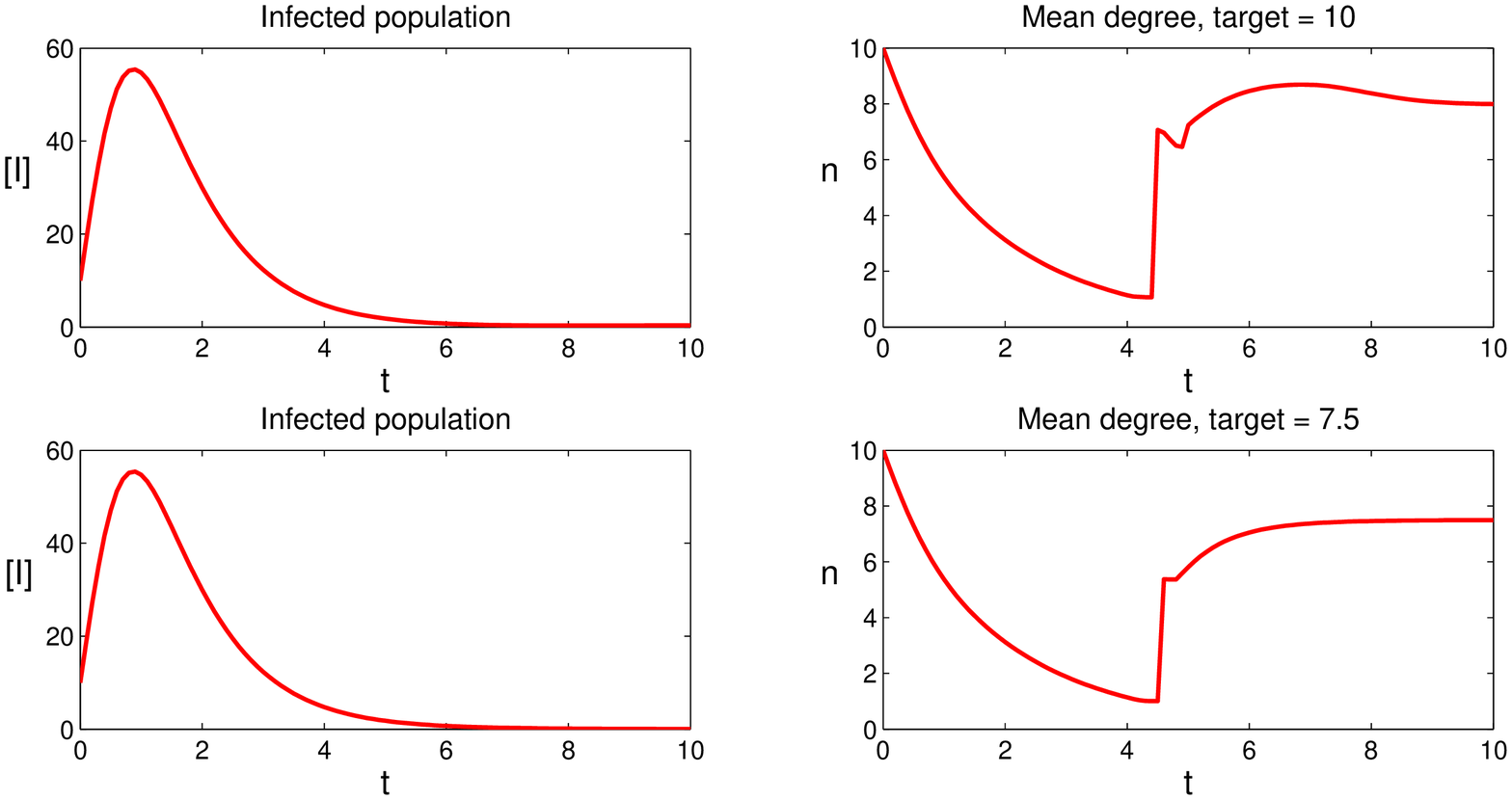}
    \caption{The effect of adjusting the control targets in terms of the time evolution of prevalence  and network connectivity for $\tau=1$, $\gamma=1$, $N=1000$, $I(0)=10$, $I^*=0$, $n(0)=10$, $M_1=6$, $M_2=0.5$, $\Delta t=0.1$, $\lambda_1=10^4$ and $\lambda_2=\lambda_3=\lambda_4=1$.} \label{fig:fig6}
\end{figure}


\begin{thebibliography}{99}
\addcontentsline{toc}{chapter}{Bibliography}

\bibitem{RMA_RMM_1991}
R.M. Anderson \& R.M. May (1991) Infectious Diseases of Humans. Oxford: \textit{Oxford University Press}.

\bibitem{JC_KAJW_KT_2013}
J. Clarke, K.A. Jane White \& Katy Turner (2013) Approximating Optimal Controls for Networks
when There Are Combinations of Population-Level and Targeted Measures Available: Chlamydia Infection as a Case-Study. \textit{Bull. Math. Biol.} \textbf{75}, 1747Ð1777.

\bibitem{RFC_HZ1995}
R.F. Curtain, H. Zwart (1995) An Introduction to Infinite-Dimensional
Linear Systems Theory, volume 21 of Texts in Applied Mathematics,
\textit{Springer-Verlag}, New York.

\bibitem{Danon}
L. Danon, A. P. Ford, T. House, C. P. Jewell, M. J. Keeling, G. O. Roberts, J. V. Ross \& M. C. Vernon (2011) "Networks and the Epidemiology of Infectious Disease", \textit{Interdisciplinary Perspectives on Infectious Diseases} \textbf{2011:284909} special issue "Network Perspectives on Infectious Disease Dynamics".

\bibitem{Grune}
L. Gr\"{u}ne, J. Pannek. Nonlinear Model Predictive Control (2011) \textit{Springer}.

\bibitem{EH_TD_11}
E. Hansen \& T. Day (2011) Optimal control of epidemics with limited resources. \textit{J. Math. Biol.} \textbf{62}, 423Ð451.


\bibitem{HouseUnifying}
T. House \& M. J. Keeling (2011) Insights from unifying modern approximations to infections on networks. \textit{J. Roy. Soc. Interface} \textbf{8}, 67-73.

\bibitem{Keeling1999}
M.J. Keeling (1999) The effects of local spatial structure on epidemiological invasions. \textit{Proc. R. Soc. Lond. B} \textbf{266}, 859-867.

\bibitem{KeelingEames}
M.J. Keeling \& K.T.D. Eames (2005) Networks and epidemic models, \textit{J. Roy. Soc. Interface} \textbf{2}, 295-307.

\bibitem{Kiss_CT}
I. Z. Kiss, D. M. Green and R. R. Kao (2005) Disease contact tracing in random and clustered networks. \textit{Proc. R. Soc. B} \textbf{272}, 1407 - 1414.

\bibitem{Barabasi}
Liu, Y. Y., Slotine, J. J. , Barab\'{a}si, A. L. (2011) Controllability of complex networks.
\textit{Nature}, \textbf{473}, 167–-173.

\bibitem{Magni}
L. Magni, D. Raimondo, F. Allg\"ower, (eds) (2009) Nonlinear Model Predictive Control -- Towards New Challenging Applications, Lecture Notes in Control and Information Sciences, \textit{Springer-Verlag}.

\bibitem{MR_WKH_74}
R. Morton \& K.H. Wickwire (1974) On the optimal control of a deterministic epidemic. \textit{Adv. Appl. Probab.}
\textbf{6}, 622Ð635.

\bibitem{EDS_1998}
E.D. Sontag (1998) Mathematical Control Theory, volume 6 of Texts in Applied
Mathematics, second edition, \textit{Springer-Verlag}, New York.

\bibitem{szab2012detailed}
A. Szab\'{o}, P. Simon \& I. Z. Kiss (2012) Detailed study of bifurcations in an epidemic model on a dynamic network.\textit{ Differ. Equ. Appl.} \textbf{4},  277-296.

\end{thebibliography}
\end{document}